\documentclass[a4paper, 12pt]{amsart}
\usepackage{amsmath, amsthm, amscd, amssymb, amsfonts, amsxtra, amssymb, latexsym}
\usepackage{enumerate}
\usepackage{verbatim}
\usepackage{tgpagella}
\usepackage{mathpazo}
\usepackage{pb-diagram}
\usepackage{hyperref}
\hypersetup{colorlinks = true,	allcolors  = blue}

\usepackage[T1]{fontenc} 
\usepackage{t1enc}
\usepackage[utf8x]{inputenc}
\usepackage{array}

\usepackage{tikz} 
\usetikzlibrary{shapes.geometric}

\theoremstyle{definition}

\theoremstyle{remark}

\usepackage{color}

\begin{document}
\numberwithin{equation}{section}
\title{A geometric proof that $\sqrt{3}$, $\sqrt{5}$ and $\sqrt{7}$ are irrational}
\author{Ricardo A.\@ Podest\'a}

\email{podesta@famaf.unc.edu.ar}

\address{Ricardo A.\@ Podest\'a -- CIEM, Universidad Nacional de C\'ordoba, CONICET, FaMAF. \linebreak  
	Av.\@ Medina Allende 2144, Ciudad Universitaria (5000) C\'ordoba, República Argentina. }

\maketitle

\begin{abstract}
We give a geometric proof that $\sqrt n$ is irrational for $n=3,5,7$ by adapting Tennenbaum's geometric proof that $\sqrt 2$ is irrational. We also show that this method cannot be used to prove the irrationality of $\sqrt n$ for a bigger $n$.
\end{abstract}

\section{Introduction}
Some of the first important mathematical results that people learn from school are the Pythagorean theorem, the fact that there are infinitely many primes, and that the number $\sqrt 2$ is irrational. They are very old theorems, dating back to the ancient Greeks, which are simultaneously deep, simple and elegant. 
Each of these statements have many, several demonstrations in the literature.
I have always found exciting the discovering of new proofs for them, 
each showing some new aspect or idea. 
   
Everyone knows the classic proof, or some variant, that $\sqrt 2$ is an irrational number. Suppose $\sqrt 2 = \frac pq$ with $p,q$ integers and $\frac pq$ a reduced fraction, that is $(p,q)=1$. Hence
\begin{equation} \label{2qp}
2q^2 =p^2
\end{equation}
and from here one deduces that $p$ is even, say $p=2m$. Then $2q^2 =4m^2$, i.e.\@ $q^2=2m^2$. 
This implies that $q$ is even also, and hence $(p,q)\ge 2$, which is a contradiction.

A surprising geometric proof of this fact was given by Tennenbaum in the 1950's (\cite{Te}). 
The proof was popularized in the 1990's by John Conway in his article 
\textit{The power of mathematics} (\cite{Co}). The proof goes as follows:
assume that $p$ and $q$ are minimal integers satisfying \eqref{2qp}. 
Interpret this geometrically as having a big square of integral size $p$ and area $p^2$ and two smaller squares of integral size $q$ and area $q^2$, \textsl{minimal} with this property, i.e.
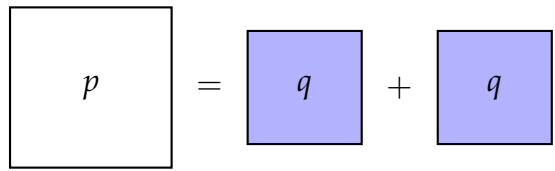
\begin{center} 
	\begin{tikzpicture}
	[scale =.125]
	\draw[thick]  (0,0) rectangle (17,17) (8.5,8.5) node {$p$} (21,8.5) node{$=$}; 
	\filldraw[fill=blue!30, thick] (25,2.5) rectangle (37,14.5) (31,8.5) node {$q$} (41,8.5) node{$+$};
	\filldraw[fill=blue!30, thick] (45,2.5) rectangle (57,14.5) (51,8.5) node {$q$};
	\end{tikzpicture} 
\end{center}
Thus, the area of the big square is the same as the area of the two smaller squares together. 
Now, placing the small squares inside the bigger one in opposite corners (see the figure below) 
they overlap forming another square of size 
$q-(p-q)=2q-p$ or alternatively $p-2(p-q)=2q-p$.

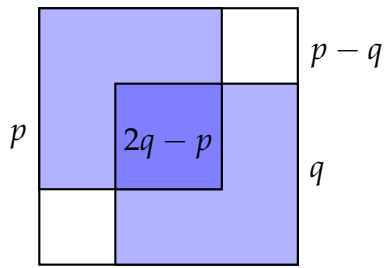
\begin{center} 
	\begin{tikzpicture}
	[scale =.2]
	\draw[thick] (0,0) rectangle (17,17) (0,8.5) node[left] {$p$}; 
	\filldraw[fill=blue!30, thick]  (5,0) rectangle (17,12) (17,14) node[right] {$p-q$}; 
	\filldraw[fill=blue!30, thick]  (0,5) rectangle (12,17) (17,6) node[right] {$q$}; 
	\filldraw[fill=blue!50, thick]  (5,5) rectangle (12,12) (8.5,8.) node {$2q-p$}; 
	\end{tikzpicture} 
\end{center}

\noindent The area of this overlapping square is counted twice and hence it must equal the area of the two not covered white little squares remaining on the corners. Hence, 
$$(2q-p)^2= 2(p-q)^2.$$ 
Graphically,
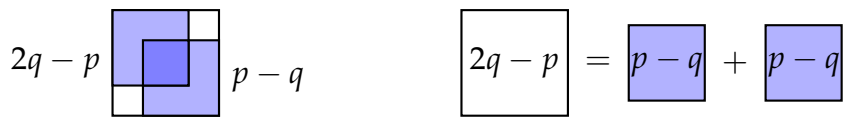
\begin{center} 
	\begin{tikzpicture}
	[scale =.2]
	\draw[thick] (0,0) rectangle (7,7) (0,3.5) node[left] {$2q-p$}; 
	\filldraw[fill=blue!30, thick]  (2,0) rectangle (7,5) (7,2.5) node[right] {$p-q$}; 
	\filldraw[fill=blue!30, thick]  (0,2) rectangle (5,7); 
	\filldraw[fill=blue!50, thick]  (2,2) rectangle (5,5); 
	\end{tikzpicture} 
%
\hspace{4em}  
\begin{tikzpicture}
[scale =.2]
\draw[thick]  (0,0) rectangle (7,7) (3.5,3.5) node {$2q-p$} (9,3.5) node{$=$}; 
\filldraw[fill=blue!30, thick] (11,1) rectangle (16,6) (13.5,3.5) node {$p-q$} (18,3.5) node{$+$};
\filldraw[fill=blue!30, thick] (20,1) rectangle (25,6) (22.5,3.5) node {$p-q$};
\end{tikzpicture} 
\end{center} 
Since $2q-p$ and $p-q$ are integers satisfying \eqref{2qp} which are smaller than $p$ and $q$, we have a contradiction. 
Thus, $\sqrt 2$ must be irrational.

In Alexander Bogomolny's celebrated blog \textit{Cut-the-knot} (\cite{Bo}) one can find more than 30 proofs, of both algebraic and geometric kind, of the irrationality of $\sqrt 2$.

\section{Geometric proof that $\sqrt{3}$ and $\sqrt 5$ are irrational}
 In a rather recent article (\cite{MM}, see also \cite{MM2}) Miller and Montague, inspired by the square-based Tennenbaum's demonstration, give beautiful geometric proofs that $\sqrt 3$ and $\sqrt 5$ are irrational numbers using overlapping areas of regular triangles and pentagons, respectively.  Next, we will prove that $\sqrt{3}$ and $\sqrt 5$ are irrational numbers by using squares a l\`a Tennenbaum.

\subsubsection*{$\sqrt 3$ is irrational}
Let us first prove the easiest case, that is that $\sqrt 3$ is irrational. This is equivalent to the fact that 
\begin{equation} \label{3qp}
3q^2 = p^2
\end{equation}
with $p>q$ integers.  
Interpret these numbers as the areas of a square of size $p$ and three squares of size $q$ 
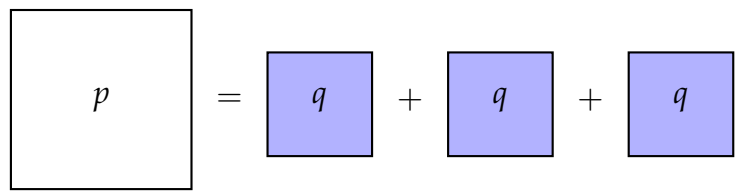
\begin{center} 
	\begin{tikzpicture}
	[scale =.125]
	\draw[thick]  (0,0) rectangle (19,19) (9.5,9.5) node {$p$} (23,9.5) node{$=$}; 
	\filldraw[fill=blue!30, thick] (27,3.5) rectangle (38,14.5) (32.5,9.5) node {$q$} (42,9.5) node{$+$};
	\filldraw[fill=blue!30, thick] (46,3.5) rectangle (57,14.5) (51.5,9.5) node {$q$} (61,9.5) node{$+$};
	\filldraw[fill=blue!30, thick] (65,3.5) rectangle (76,14.5) (70.5,9.5) node {$q$};
	\end{tikzpicture} 
\end{center}
and assume they are minimal with this property.

Inside the bigger square we place the three smaller squares sitting on a diagonal, two on the opposite corners and the remaining one at the center: 

\

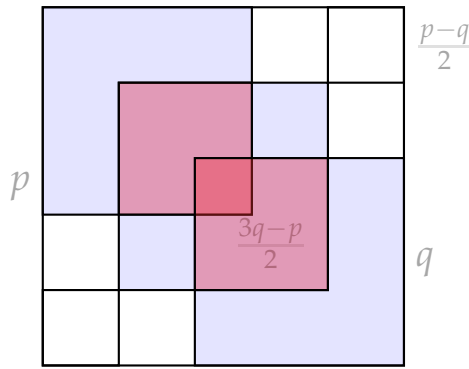
\begin{center} 
	\begin{tikzpicture}
	[scale =.25]
	\pgfsetfillopacity{0.35}
	\draw[thick] (0,0) rectangle (19,19) (0,9.5) node[left] {\large {$p$}}; 
	\draw[thick] (0,0) rectangle (4,4);
	\draw[thick] (15,15) rectangle (19,19) (19,17) node[right] {\large $\frac{p-q}2$};
	\filldraw[fill=blue!30, thick]  (8,0) rectangle (19,11) (19,5.5) node[right] {\large $q$}; 
	\filldraw[fill=blue!30, thick]  (0,8) rectangle (11,19) (17,6) ; 
	\filldraw[fill=blue!30, thick]  (4,4) rectangle (15,15) (17,6); 
	\filldraw[fill=red!80, thick]   (4,8) rectangle (11,15) (12,6.5) node {\large $\frac{3q-p}2$}; 
	\filldraw[fill=red!80, thick]   (8,4) rectangle (15,11) (8.5,8.); 
	\end{tikzpicture} 
\end{center}

\noindent 
Due to the overlapping, we now have squares of some different sizes. 
The overlapped area (in red), which is counted twice, is the same as the area of the white squares not touched by the three original squares (in blue). There are 6 white squares of size $\frac{p-q}2$ 
on the corners, whose total area equals the area of the two red squares on the diagonal of size 
$$q-\tfrac{p-q}2 = \tfrac{3q-p}2.$$    
In this way we have $6(\frac{p-q}2)^2 = 2 (\tfrac{3q-p}2)^2$, that is 
$$3(\tfrac{p-q}2)^2 = (\tfrac{3q-p}2)^2.$$
By \eqref{3qp} the parity of $p$ and $q$ are the same. 
Thus, we have found smaller squares of integral sizes $\frac{p-q}2$ and $\tfrac{3q-p}2$ satisfying \eqref{3qp}, which is absurd.

\subsubsection*{$\sqrt 5$ is irrational}
We now move to the next case $\sqrt 5$. Suppose there are integers $p$ and $q$ such that 
\begin{equation} \label{5qp}
5q^2=p^2.
\end{equation}
That is, there is a square of integral size $p$ and area $p^2$ which equal the total area of 5 squares of integral size $q$ and individual area $q^2$. 
Place the five squares of size $q$ symmetrically equidistributed along the principal diagonal: 

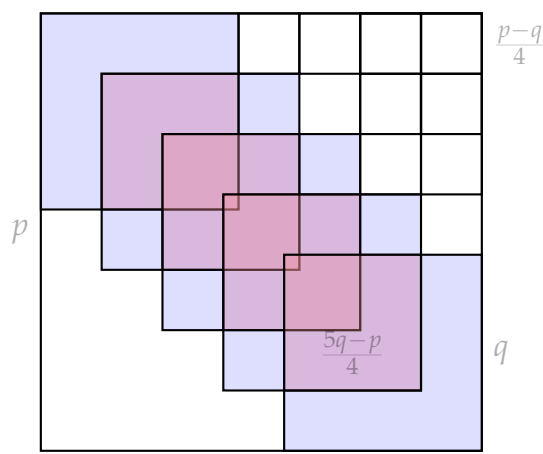
\begin{center} 
	\begin{tikzpicture}
	[scale =.2]
	\pgfsetfillopacity{0.325}
	\draw[thick] (0,0) rectangle (29,29) (0,14.5) node[left] {$p$}; 
	\draw[thick] (25,25) rectangle (29,29) (29,27) node[right] {$\frac{p-q}4$};
	\draw[thick] (21,25) rectangle (25,29);
	\draw[thick] (17,25) rectangle (21,29);
	\draw[thick] (13,25) rectangle (17,29);
	\draw[thick] (21,21) rectangle (25,25);	
	\draw[thick] (25,21) rectangle (29,25);
	\draw[thick] (21,17) rectangle (25,29);
	\draw[thick] (25,17) rectangle (29,29);
	\filldraw[fill=blue!40, thick]  (16,0) rectangle (29,13) (29,6.5) node[right] {$q$}; 
	\filldraw[fill=blue!40, thick]  (0,16) rectangle (13,29);  
	\filldraw[fill=blue!40, thick]  (4,12) rectangle (17,25); 
	\filldraw[fill=blue!40, thick]  (8,8)  rectangle (21,21); 
	\filldraw[fill=blue!40, thick]  (12,4) rectangle (25,17);
	\filldraw[fill=red!40, thick]  (4,16) rectangle (13,25); 
	\filldraw[fill=red!40, thick]  (8,12) rectangle (17,21);
	\filldraw[fill=red!40, thick]  (12,8) rectangle (21,17);	
	\filldraw[fill=red!40, thick]  (16,4) rectangle (25,13) (20.5,6.5) node{$\frac{5q-p}4$};	
	\end{tikzpicture} 
\end{center}

On the upper-right corner there are $10$ squares of size $\tfrac{p-q}4$ (hence a total of $20$ such white squares). 
There are 4 overlapping squares (in red) of size 
$$q-\tfrac{p-q}4=\tfrac{5q-p}4.$$
The area of each one of the 4 red squares is counted twice.   
Since the area of the overlapping squares must equal the area of the white squares, we have 
$20 (\tfrac{p-q}4)^2= 4 (\tfrac{5q-p}4)^2$, that is
\begin{equation} \label{5qpred}
5(\tfrac{p-q}4)^2= (\tfrac{5q-p}4)^2.
\end{equation}
If $\tfrac{p-q}4$ and $\tfrac{5q-p}4$ are integers then we have constructed squares of integral size satisfying \eqref{5qp} which are smaller than the original ones, contradicting the assumption. On the other hand, if $\tfrac{p-q}4$ and $\tfrac{5q-p}4$ are not integers, \eqref{5qpred} is equivalent to 
\begin{equation*} \label{5qpred2}
5(\tfrac{p-q}{2})^2= (\tfrac{5q-p}{2})^2
\end{equation*}
and the same conclusion as before can be deduced from here. That is, the squares in the drawing are not integral, but we can find integral squares smaller than the original ones, satisfying \eqref{5qp}. It is clear that $\frac{p-q}2$ and $\frac{5q-p}2$ are integers since $p$ and $q$ have the same parity by \eqref{5qp}. Also, 
$$\tfrac{p-q}2 < p-q < p \qquad \text{and} \qquad \tfrac{5q-p}2<p$$ since 
this is equivalent to $5q < 3p = 3\sqrt 5 q$ which holds since $5<3\sqrt 5$.

\section{Can this method be generalized to $\sqrt{n}$?}
It may seem that we can generalize the geometric method of Tennenbaum previously used to show that $\sqrt n$ is irrational for any $n$ not a square. Let us see what happens.

Let $n$ be a non-square and proceed as before. Suppose that $\sqrt n$ is rational. That is, there are integers $p,q$ such that  
\begin{equation} \label{nqp}
nq^2=p^2.
\end{equation}
In other words, there is a square of integral size $p$ with the same area of $n$ integral squares of size $q$.
Equidistribute the $n$ smaller squares with their centers on the diagonal of the square of size $p$. 

Note that the area of the big square not covered by the smaller squares can be divided into little squares of size
$$ \tfrac{p-q}{n-1}$$ 
and that the number of these little squares on the upper-right corner is exactly 
\begin{equation} \label{tn}
t_{n-1} = 1 + 2 + \cdots +{n-1} = \tfrac{n(n-1)}2 = \tbinom{n}{2}
\end{equation}
where $t_n$ is the $n$-th triangular number. In the previous cases of Section 2 we have 
$t_2=1+2=3$ for $n=3$ and $t_4=1+2+3+4=10$ for $n=5$. 
On the other hand, there are $n-1$ overlapping squares each of size 
$$q-\tfrac{p-q}{n-1} = \tfrac{nq-p}{n-1}.$$
Taking into account the equality of the areas we have 
$$ 2 t_{n-1} \big(\tfrac{p-q}{n-1}\big)^2 = (n-1) \big( \tfrac{nq-p}{n-1} \big)^2$$
which, by \eqref{tn}, is equivalent to 
$2 \tfrac{(n-1)n}{2} (p-q)^2 = (n-1)(nq-p)^2,$ 
that is
\begin{equation} \label{nqpred}
n (p-q)^2 = (nq-p)^2.
\end{equation}
Is this an absurd? Can we deduce from here the irrationality of $\sqrt n$? 
Well, to get a contradiction we need that $nq-p<p$, but this is equivalent to $nq <2p = 2\sqrt n q$ which holds if and only if $n<2\sqrt n$. Equivalently, $n^2<4n$ which only holds for $n=2,3$.

Fortunately, we still have the trick used for $\sqrt 5$. Since $p$ and $q$ have the same parity, if we consider $n$ to be odd we have that both $\frac{p-q}2$ and $\frac{nq-p}2$ are integers. Also,  
equation \eqref{nqpred} is equivalent to 
\begin{equation} \label{nqpred2}
n (\tfrac{p-q}2)^2 = (\tfrac{nq-p}2)^2.
\end{equation}
Thus, we get a contradiction if we prove that $\frac{nq-p}2 < p$. But this is equivalent to 
$$nq < 3p = 3\sqrt n q \quad \Leftrightarrow \quad n<3 \sqrt n \quad \Leftrightarrow \quad n(n-9)<0$$
which holds for $n<9$, that is $n=3,5,7$. In this way, we have also proved that $\sqrt 7$ is irrational with Tennenbaum's method.
Can we extend the trick used to get further values of $n$? 
No! Because, if we take for instance $\frac{p-q}4$ and $\frac{nq-p}4$ there is no way to ensure that these numbers are integers. 


\section{Final remarks}
For the sake of completeness, we mention that there are some few geometrical proofs that $\sqrt n$ is irrational for every $n$ not a square.
We point out the proof of Derek
Ball based on successively cutting squares from rectangles (\cite{DB}), the one of Terence Jackson using similar triangles (\cite{TJ}) and a recent one from Nick Lord tiling rectangles (\cite{NL}). 

%
%

\end{document}